\documentclass[a4paper,11pt,english]{article}
\usepackage[utf8]{inputenc}
\usepackage{babel,amsmath,amssymb,amsthm,enumitem,graphicx}
\usepackage[sort&compress,numbers,square]{natbib}
\usepackage[colorlinks,allcolors=blue]{hyperref}
\usepackage[margin=1.25in]{geometry}
\usepackage[small,bf]{titlesec}
\titlelabel{\thetitle.\hspace{0.5em}}

\newtheorem*{theorem}{Main Theorem}
\newtheorem*{lemma}{Lemma}
\theoremstyle{definition}
\newtheorem*{remark}{Remark}
\newtheorem*{example}{Example}

\DeclareMathOperator{\E}{E}
\DeclareMathOperator{\sgn}{sgn}

\DeclareMathOperator{\I}{I}
\renewcommand{\P}{\mathrm{P}}
\newcommand{\F}{\mathcal{F}}
\newcommand{\FF}{\mathbb{F}^X}

\renewcommand{\tilde}{\widetilde}
\renewcommand{\epsilon}{\varepsilon}
\newcommand{\prt}[2]{\frac{\partial #1}{\partial #2}}

\title{A sequential test for the drift of a Brownian motion\\
with a possibility to change a decision}
\author{Mikhail Zhitlukhin\thanks{Steklov Mathematical Institute of the
Russian Academy of Sciences. 8 Gubkina St., Moscow, Russia.
Email: mikhailzh@mi-ras.ru.
The research was supported by the Russian Science Foundation,
project no.\ 19-11-00290.}}
\date{25 July 2020}

\begin{document}
\maketitle

\begin{abstract}
We construct a Bayesian sequential test of two simple hypotheses about the
value of the unobservable drift coefficient of a Brownian motion, with a
possibility to change the initial decision at subsequent moments of time for
some penalty. Such a testing procedure allows to correct the initial
decision if it turns out to be wrong. The test is based on observation of
the posterior mean process and makes the initial decision and, possibly,
changes it later, when this process crosses certain thresholds. The solution
of the problem is obtained by reducing it to joint optimal stopping and
optimal switching problems.


\end{abstract}

\section{Introduction}
We consider a problem of sequential testing of two simple hypotheses about
the value of the unknown drift coefficient of a Brownian motion. In usual
sequential testing problems (see e.g.\ the seminal works
\cite{WW48,S67,IS84} or the recent monographs \cite{BLS12,TNB14}), a testing
procedure must be terminated at some stopping time and a decision about the
hypotheses must be made. In contrast, in the present paper we propose a new
setting, where a testing procedure does not terminate and it is allowed to
change the initial decision (for the price of paying some penalty) if, given
later observations, it turns out that it is incorrect.

We will work in a Bayesian setting and assume that the drift coefficient has
a known prior distribution on a set of two values. A decision rule consists
of an initial decision $(\tau,d)$, where $\tau$ is the moment at which the
decision is made and $d$ is a two-valued function showing which hypothesis
is accepted initially, and a sequence of stopping times $\tau_n$, at which
the decision can be changed later. The goal is to minimize a penalty
function which consists of the three parts: a penalty for the waiting time
until the initial decision, a penalty for a wrong decision proportional to
the time during which the corresponding wrong hypothesis is being accepted,
and a penalty for each change of a decision.

This study was motivated by the paper \cite{MUZ19}, where a sequential
multiple changepoint detection problem was considered. That problem consists
in tracking of the value of the unobservable drift coefficient of a Brownian
motion, which is modeled by a telegraph process (a two-state Markov process)
switching between $-1$ and $+1$ at random times. In the present paper, we
deal with a similar tracking procedure and a penalty function, but the
difference is that the unobservable drift coefficient does not change. Among
other results on multiple changepoint detection, one can mention the paper
\cite{G15}, where a tracking problem for a general two-state Markov process
with a Brownian noise was considered, and the paper \cite{BL09}, which
studied a tracking problem for a compound Poisson process.

We solve our problem by first representing it as a combination of an optimal
stopping problem and an optimal switching problem (an optimal switching
problem is an optimal control problem where the control process assumes only
two values). The optimal stopping problem allows to find the initial
stopping time, while the subsequent moments when the decision is changed are
found from the optimal switching problem. Consequently, the value function
of the optimal switching problem becomes the payoff function of the optimal
stopping problem. Then both of the problems are solved by reducing them to
free-boundary problems associated with the generator of the posterior mean
process of the drift coefficient. We consider only the symmetric case (i.e.\
type I and type~II errors are of the same importance), in which the solution
turns out to be of the following structure. First an observer waits until
the posterior mean process exists from some interval $(-A,\,A)$ and at that
moment of time makes the initial decision. Future changes of the decision
occur when the posterior mean process crosses some thresholds $-B$ and $B$.
The constants $A,B$ are found as unique solutions of certain equations.

The rest of the paper consists of the three sections:
Section~\ref{section-problem} describes the problem,
Section~\ref{section-result} states the main theorem which provides the
optimal decision rule, Section~\ref{section-proof} contains its proof.

\section{The model and the optimality criterion}
\label{section-problem}

Let $(\Omega,\F, \P)$ be a complete probability space. Suppose one can
observe a process $X_t$ defined on this probability space by the relation
\begin{equation}
X_t = \mu \theta t + B_t,\label{X}
\end{equation}
where $B_t$ is a standard Brownian motion, $\mu>0$ is a known constant, and
$\theta$ is a $\pm1$-valued random variable independent of $B_t$. It is
assumed that neither $\theta$ nor $B_t$ can be observed directly. The goal
is to find out whether $\theta=1$ or $\theta=-1$ by observing the process
$X_t$ sequentially. Note that the case when the drift coefficient of $X_t$
can take on two arbitrary values $\mu_1\neq \mu_2$ can be reduced
to~\eqref{X} by considering the process $X_t - \frac12(\mu_1+\mu_2)t$.

We will assume that the prior distribution of $\theta$ is known and is
characterized by the probability $p = \P(\theta=1)$. Recall that usual
settings of sequential testing problems consist in that an observer must
choose a stopping time $\tau$ of the (completed and right-continuous)
filtration $\FF = (\F_t^X)_{t\ge 0}$ generated by $X_t$, at which the
observation is stopped, and an $\F_\tau^X$-measurable function $d$ with
values $-1$ or $+1$ that shows which of the two hypotheses is accepted at
time $\tau$. The choice of $(\tau,d)$ depends on a particular optimality
criterion which combines penalties for type I and type II errors, and a
penalty for observation duration. But, in any case, a test terminates at
time $\tau$.

In this paper we will focus on a setting where an observer can change a
decision made initially at time $\tau$ and the testing procedure does not
terminate.

By a \emph{decision rule} we will call a triple $\delta=(\tau_0,d,T)$, where
$\tau_0$ is an $\FF$-stopping time, $d$ is an $\F_{\tau_0}^X$-measurable
function which assumes values $\pm1$, and $T=(\tau_1,\tau_2\ldots)$ is a
sequence of $\FF$-stopping times such that $\tau_n\le\tau_{n+1}$ for all
$n\ge 0$. At the moment $\tau_0$, the initial decision $d$ is made. Later,
if necessary, an observer can change the decision to the opposite one, and
the moments of change are represented by the sequence $T$. Thus, if, for
example, $d=1$, then at $\tau_0$ an observer decides that $\theta=1$ and at
$\tau_1$ switches the opinion to $\theta=-1$; at $\tau_2$ switches back to
$\theta=1$, and so on. It~may be the case that $\tau_n=+\infty$ starting
from some $n$; then the decision is changed only a finite number of times
(the optimal rule we construct below will have this property with
probability 1).

With a given decision rule $\delta$, associate the $\FF$-adapted process
$D_t^\delta$ which expresses the current decision at time $t$,
\[
D_t^\delta =
\begin{cases}
0, &\text{if}\ t< \tau_0,\\
d, &\text{if}\ t\in [\tau_{2n}, \tau_{2n+1}),\\
-d, &\text{if}\ t\in [\tau_{2n+1}, \tau_{2n+2}),
\end{cases}
\]
and define the \emph{Bayesian risk function} 
\begin{equation}
R(\delta) = \E \biggl(c_0 \tau_0 + c_1 \int_{\tau_0}^\infty \I(D_t^\delta
\neq \theta) dt + c_2 \sum_{t> \tau_0} \I(D_{t-}^\delta \neq
D_t^\delta)\biggr),\label{R}
\end{equation}
where $c_i>0$ are given constants.

The problem that we consider consists in finding a decision rule $\delta^*$
which minimizes $R$, i.e.
\[
R(\delta^*) = \inf_{\delta} R(\delta).
\]
Such a decision rule $\delta^*$ will be called \emph{optimal}.

One can give the following interpretation to the terms under the expectation
in \eqref{R}. The term $c_0\tau_0$ is a penalty for a delay until making the
initial decision. The next term is a penalty for making a wrong decision,
which is proportional to the time during which the wrong hypothesis is being
accepted. The last term is a penalty for changing a decision, in the amount
$c_2$ for each change. Note that the problem we consider is symmetric (i.e.\
type I and type II errors are penalized in the same way); in principle, an
asymmetric setting can be studied as well.

\section{The main result}
\label{section-result}
To state the main result about the optimal decision rule, introduce the
posterior mean process
\[
M_t = \E(\theta \mid \F_t^X).
\]
As follows from known results, the process $M_t$ satisfies the stochastic
differential equation
\begin{equation}
d M_t = \mu(1-M_t^2) d \tilde B_t, \qquad M_0=2p-1,\label{M}
\end{equation}
where $\tilde B_t$ is a Brownian motion with respect to $\FF$ (an
\emph{innovation process}, see, e.g., Chapter~7~in~\cite{LS01}), which
satisfies the equation
\[
d \tilde B_t = d X_t - M_t d t.
\]
Representation~\eqref{M} can be obtained either directly from filtering
theorems (see Theorem~9.1~in~\cite{LS01}) or from the known equation for the
posterior probability process $\pi_t = \P(\theta=1\mid \F_t^X)$ (see
Chapter~VI~in~\cite{PS06}) since $M_t = 2\pi_t - 1$. In the explicit form,
$M_t$ can be expressed through the observable process $X_t$ as
\[
M_t = 1 - \frac{2(1-p)}{{p} e^{2\mu X_t}+1-p}.
\]

Introduce the two thresholds $A,B \in(0,1)$, which depend on the parameters
$\mu,c_0,\allowbreak c_1,c_2$ of the problem, and will define the switching
boundaries for the optimal decision rule. The threshold $B$ is defined as
the solution of the equation
\begin{equation}
\ln\frac{1-B}{1+B} + \frac{2B}{1-B^2} = \frac{2\mu^2c_2}{c_1},\label{B}
\end{equation}
and the threshold $A$ is defined as the solution of the equation
\begin{equation}
\biggl(\frac{c_1}{2c_0} -1\biggr) \ln\frac{1-A}{1+A} + \frac{2}{1+A}
\biggl(\frac{c_1}{2c_0} + \frac{A}{1-A} \biggr) =
\frac{c_1}{c_0(1-B^2)}.\label{A}
\end{equation}
The next simple lemma shows that $A$ and $B$ are well-defined. Its proof is
rather straightforward and is omitted.
\begin{lemma}
\label{lemma}
Equations \eqref{B},~\eqref{A} have unique solutions $A,B\in(0,1)$. If
$c_1=2c_0$, then $A=B$.
\end{lemma}

The following theorem, being the main result of the paper, provides the
optimal decision rule in an explicit form.
\begin{theorem}
\label{theorem}
The optimal decision rule $\delta^*=(\tau_0^*,d^*,T^*)$ consists of the
stopping time $\tau_0^*$ and the decision function $d^*$ defined by the
formulas
\[
\tau_0^* = \inf\{t\ge 0: |M_t| \ge A\}, \qquad d^* = \sgn M_{\tau_0},
\]
and the sequence of stopping times $T^*=(\tau_n^*)_{n= 1}^\infty$ which on
the event $\{d^*=1\}$ are defined by the formulas
\begin{equation}
\tau_{2k+1}^* = \inf\{t\ge \tau_{2k}^*: M_t \le -B\}, \quad \tau_{2k+2}^* =
\inf\{t\ge \tau_{2k+1}^*: M_t \ge B\},\label{tau-1}
\end{equation}
and on the event $\{d=-1\}$ by the formulas
\begin{equation}
\tau_{2k+1}^* = \inf\{t\ge \tau_{2k}^*: M_t \ge B\}, \quad
\tau_{2k+2}^* = \inf\{t\ge \tau_{2k+1}^*: M_t \le -B\} \label{tau-2}
\end{equation}
(where $\inf\emptyset = +\infty$).
\end{theorem}

\begin{example}
Figure~\ref{figure-process} illustrates how the optimal decision rule works.
In this example, we take $p=0.5$, $\mu=1/3$, $c_0=2/3$, $c_1=1$, $c_2=3/2$.
The thresholds $A,B$ can be found numerically, $A\approx 0.37$,
$B\approx0.55$.

The simulated path on the left graph has $\theta=1$. The rule $\delta^*$
first waits until the process $M_t$ exists from the interval $(-A,A)$. Since
in this example it exists through the lower boundary (at $\tau_0^*$), the
initial decision is $d^*=-1$ (incorrect). Then the rule waits until $M_t$
crosses the threshold $B$, and changes the decision to $\theta=1$ at
$\tau_1^*$.

\begin{figure}[h]
\centering \includegraphics[width=14cm]{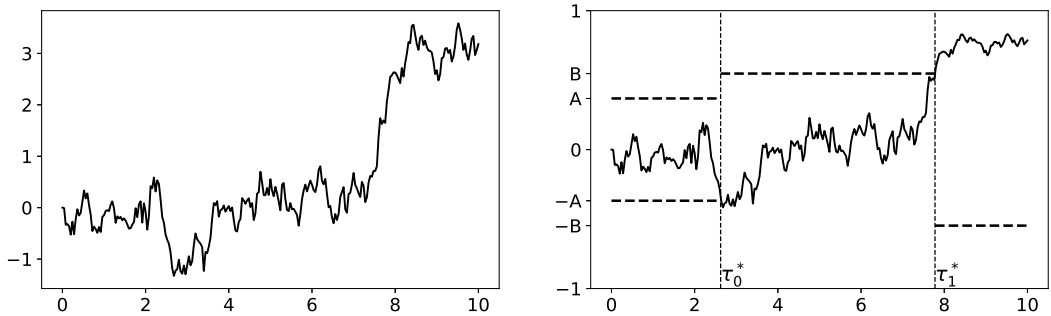}
\caption{Left: the process $X_t$; right: the process $M_t$. Parameters:
$p=0.5$, $\mu=1/3$, $c_0=2/3$, $c_1=1$, $c_2=3/2$.}
\label{figure-process}
\end{figure}
\end{example}

\section{Proof of the Main Theorem}
\label{section-proof}

Let us denote by $\P_x$ and $\E_x$ the probability measure and the
expectation under the assumption $\P(\theta=1)=(x+1)/2$, so the posterior
mean process $M_t$ starts from the value $M_0=x$. It is easy to verify that
\[
\P_x ( D_t^\delta \neq \theta \mid \F_{t}^X) = \frac{1- M_tD_t^\delta}{2},
\]
and, by taking intermediate conditioning with respect to $\F_t^X$ in
\eqref{R}, we can see that we need to solve the problem
\begin{equation}
V^*(x) = \inf_\delta \E_x \biggl(c_0 \tau_0 + \frac{c_1}{2}
\int_{\tau_0}^\infty (1-M_tD_t^\delta) dt + c_2 \sum_{t> \tau_0}
\I(D_{t-}^\delta \neq D_t^\delta)\biggr), \qquad x\in[-1,1]
\label{Vstar}
\end{equation}
(by ``to solve'' we mean to find $\delta$ at which the infimum is attained
for a given $x$; in passing we will also find the function $V^*(x)$ in an
explicit form).

Observe that there exists the limit $M_\infty := \lim_{t\to\infty }M_t =
\theta$ a.s. Hence the solution of problem~\eqref{Vstar} should be looked
for only among decision rules $\delta$ such that $D^\delta_t$ has a finite
number of jumps and $D_\infty^\delta = \theta$ (note that the rule
$\delta^*$ satisfies these conditions). In view of this, for a stopping time
$\tau_0$ denote by $\mathcal{D}(\tau_0)$ the class of all $\FF$-adapted
c\`adl\`ag processes $D_t$ such that, with probability 1, they assume values
$\pm 1$ after $\tau_0$, have a finite number of jumps, and satisfy the
condition $D_\infty = \theta$. Let $U^*(\tau_0)$ be the value of the
following optimal switching problem:
\begin{equation}
U^*(\tau_0) = \inf_{D\in \mathcal{D}(\tau_0)} \E_x \biggl(\frac{c_1}{2}
\int_{\tau_0}^\infty (1-M_tD_t) dt + c_2 \sum_{t> \tau_0} \I(D_{t-}
\neq D_t)\biggr).\label{Ustar}
\end{equation}
Consequently, problem \eqref{Vstar} can be written in the form
\begin{equation}
V^*(x) = \inf_{\tau_0} \E_x (c_0\tau_0 + U^*(\tau_0)).
\label{Vstar-2}
\end{equation}
Thus, to show that the decision rule $\delta^*$ is optimal, it will be
enough to show that $\tau_0^*$ delivers the infimum in the problem $V^*$,
and $D^{\delta^*}$ delivers the infimum in the problem $U^*(\tau_0^*)$. In
order to do that, we are going to use a usual approach based on ``guessing''
a solution and then verifying it using It\^o's formula. Since this approach
does not show how to actually find the functions $V^*$ and $U^*$, in the
remark after the proof we provide heuristic arguments that can be used for
that.

\begin{figure}[b]
\centering \includegraphics[width=12cm]{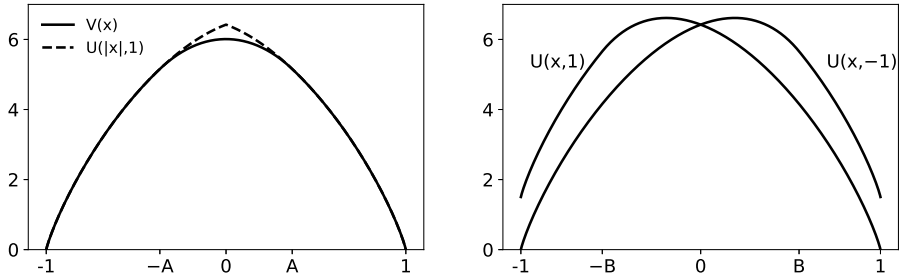}
\caption{The functions $V(x)$ and $U(x,y)$. The parameters $\mu,c_0,c_1,c_2$
 are the same as in Figure~\ref{figure-process}.}
\label{figure-VU}
\end{figure}

We will first deal with $U^*$. Let $B$ be the constant from \eqref{B}.
Introduce the ``candidate'' function $U(x,y)$ $x\in[-1,1]$, $y\in\{-1,1\}$,
defined by
\begin{alignat}{2}
&U(x,1) = \frac{c_1(1-x)}{4\mu^2}\biggl(\ln\frac{1+x}{1-x} +
\frac{2}{1-B^2}\biggr),\qquad && x\in(-B,\;1]\label{U-eq-1},\\
&U(x,1) = U(-x,1)+c_2\vphantom{\Big|},&& x\in[-1,\;-B]\label{U-eq-2},\\
&U(x,-1) = U(-x,1)\vphantom{\Big|},&& x\in [-1,\;1]\label{U-eq-3}
\end{alignat}
(see Figure~\ref{figure-VU}, which depicts the function $U(x,y)$, as well as
the function $V(x)$ defined below, with the same parameters as in the
example in the previous section).

We are going to show that $U^*(\tau_0) = U(|M_{\tau_0}|, 1)$. Let $L f$
denotes application of the generator of the process $M_t$ to a sufficiently
smooth function~$f$, i.e.
\[
L f(x)  = \frac{\mu^2}{2} (1-x^2)^2 \frac{\partial^2}{\partial x^2} f(x).
\]
By $U'$ and $\Delta U$ denote, respectively, the derivative with respect to
the first argument, and the difference with respect to the second argument
of $U$, i.e.
\[
U'(x,y) = \prt Ux(x,y), \qquad \Delta U(x,y)= U(x,y) - U(x,-y). 
\]
From the above explicit construction \eqref{U-eq-1}--\eqref{U-eq-3}, it is not
difficult to check that $U(x,y)$ has the following properties:
\begin{enumerate}[leftmargin=*,topsep=0.5em,label=(U.\arabic*)]
\item\label{U1} $U(x,y) \in C^1$ in $x$ for $x\in(-1,1)$, and $U(x,y) \in
C^2$ in $x$ except at points $x=-y B$;

\item\label{U2} $(1-x^2) U'(x,y)$ is bounded for $x\in(-1,1)$;

\item\label{U3} $LU(x,y) = -{c_1}(1-xy)/2$ if $xy > -B$, and $LU(x,y) \ge
-{c_1}(1-xy)/2$ if $xy < -B$;

\item\label{U4} $\Delta U(x,y) = -c_2$ if $xy\ge B$, and $\Delta U(x,y) \ge
-c_2$ if $xy<B$.
\end{enumerate}

Consider any process $D\in\mathcal{D}(\tau_0)$ and let $(\tau_n)_{n\ge 1}$
be the sequence of the moments of its jumps after $\tau_0$.
Property~\ref{U1} allows to apply It\^o's formula to the process $U(M_t,
D_t)$, from which for any $s>0$ we obtain
\begin{equation}
\begin{split}
U(M_{s\vee \tau_0}, D_{s\vee\tau_0}) &= U(M_{\tau_0}, D_{\tau_0}) +
\sum_{n\,:\,\tau_{n-1} \le s}\biggl( \int_{\tau_{n-1}}^{s\wedge \tau_n} L
U(M_t, D_t) \I(M_t\neq -D_t B) dt \\&+ \mu\int_{\tau_{n-1}}^{s\wedge \tau_n}
(1-M_t^2)U'(M_t,D_t) d \tilde B_t + \Delta U(M_{\tau_n},
D_{\tau_n})\I(s\ge\tau_n) \biggr).
\end{split}
\label{U-ito}
\end{equation}
Take the expectation $\E_x(\,\cdot\mid \F_{\tau_0}^X)$ of the both sides of
\eqref{U-ito}. By~\ref{U2}, the integrand in the stochastic integral is
uniformly bounded, so its expectation is zero. Passing to the limit
$s\to\infty$ and using the equality $D_\infty = M_\infty$, which implies
$U(M_{s\vee\tau_0},D_{s\vee\tau_0})\to0$ as $s\to\infty$, we obtain
\begin{equation}
U(M_{\tau_0}, D_{\tau_0}) \le
\E_x\biggl(\frac{c_1}{2}\int_{\tau_0}^\infty (1- M_tD_t) dt +
c_2\sum_{t> \tau_0}  \I( D_t \neq D_{t-}) \;\Big|\; \F_{\tau_0}^X\biggr),
\label{U-ito-2}
\end{equation}
where to get the inequality we used property~\ref{U3} for the first term
under the expectation and \ref{U4} for the second term. Taking the infimum
of the both sides of \eqref{U-ito-2} over $D \in \mathcal{D}(\tau_0)$ we
find
\begin{equation}
U(M_{\tau_0}, D_{\tau_0}) \le U^*(\tau_0).\label{U-ineq}
\end{equation}
On the other hand, if the process $D_t$ is such that $D_{\tau_0} = \sgn
M_{\tau_0}$ (let $\sgn 0 = 1$, if necessary) and its jumps after $\tau_0$
are identified with the sequence $(\tau_n)_{n\ge 1}$ defined as in
\eqref{tau-1}--\eqref{tau-2} but with arbitrary $\tau_0$ in place of
$\tau_0^*$, then we would have the equality in~\eqref{U-ito-2}, as follows
from \ref{U3} and \ref{U4}. Together with \eqref{U-ineq}, this implies that
$U^*(\tau_0) = U(M_{\tau_0}, \sgn M_{\tau_0}) = U(|M_{\tau_0}|, 1)$ and the
infimum in the definition of $U^*(\tau_0)$ is attained at this process
$D_t$.

Let us now consider the problem $V^*$. As follows from the above arguments,
we can write it in the form
\begin{equation}
V^*(x) = \inf_{\tau_0} \E_x (c_0\tau_0 + U(|M_{\tau_0}|,
1)).\label{Vstar-3}
\end{equation}
It is clear that it is enough to take the infimum only over stopping times
with finite expectation.

Let $A$ be the constant defined in \eqref{A}, and put
\begin{equation}
K = \biggl(\frac{c_1(1-A)}{4\mu^2} + \frac{c_0A}{2\mu^2}\biggr)
\ln\frac{1+A}{1-A} + \frac{c_1(1-A)}{2\mu^2(1-B^2)}.\label{K}
\end{equation}
Introduce the ``candidate'' function $V(x)$,
$x\in[-1,1]$:
\begin{alignat}{2}
&V(x) = \frac{c_0x}{2\mu^2} \ln \frac{1-x}{1+x} + K,\qquad &&|x|< A,\label{V-eq-1}\\
&V(x) = U(|x|,1),\qquad &&|x|\ge A. \label{V-eq-2}
\end{alignat}
It is straightforward to check that $V(x)$ has the following properties:
\begin{enumerate}[leftmargin=*,topsep=0.5em,label=(V.\arabic*)]
\item\label{V1}  $V(x) \in C^1$ in $x$ for $x\in(-1,1)$, and $V(x) \in C^2$
in $x$ except at points $x=\pm A$;

\item\label{V2} $(1-x^2) V'(x)$ is bounded for $x\in(-1,1)$;

\item\label{V3} $LV(x) = -c_0$ if $|x| < A$, and $LV(x) \ge
-c_0$ if $|x| > A$;

\item\label{V4} $V(x) = U(|x|,1)$ if $|x|\ge A$, and $V(x) \le U(|x|,1)$ if
$|x|<A$.
\end{enumerate}
Applying I\^o's formula to the process $V(M_t)$ and taking the expectation,
for any stopping time $\tau_0$ with $\E \tau_0<\infty$ we obtain
\[
\E_x V(M_{\tau_0}) = V(x) + \E_x \int_0^{\tau_0} LV(M_s) ds
\]
(It\^o's formula can be applied in view of \ref{V1}; the expectation of the
stochastic integral, which appears in it, is zero in view of \ref{V2} and
the finiteness of $\E \tau_0$).

From \ref{V3} and \ref{V4}, we find
\begin{equation}
V(x) \le \E_x (c_0\tau_0 + U(|M_{\tau_0}|, 1)),\label{V-ineq}
\end{equation}
so, after taking the infimum over $\tau_0$, we get $V(x) \le V^*(x)$. On the
other hand, for the stopping time $\tau_0^*$ we have the equality
in~\eqref{V-ineq}, so $V(x) = V^*(x)$. Consequently, $\tau_0^*$ solves the
problem $V^*$.

The proof is complete.

\begin{remark}
The above proof does not explain how to find the functions $V(x)$ and
$U(x,y)$. Here we provide arguments which are based on well-known ideas from
the optimal stopping theory and allow to do that. The reader is referred,
e.g., to the monograph \cite{PS06} for details.

Since the process $M_t$ is Markov, we can expect that the optimal process
$D_t$ for $U^*$ should depend only on current values of $M_t$ and $D_{t-}$.
Moreover, it is natural to assume that $D_t$ should switch from $1$ to $-1$
when $M_t$ becomes close to $-1$, and switch from $-1$ to $1$ when $M_t$
becomes close to $1$. The symmetry of the problem suggests that there should
be a threshold $B$ such that the switching occurs when $M_t$ crosses the
levels $\pm B$. This means that the optimal sequence of stopping times $T^*$
is of the form \eqref{tau-1}--\eqref{tau-2}. Consequently, in the set
$\{(x,y) : x>-yB\}$, where $x$ corresponds to the value of $M_t$ and $y$
corresponds to the value of $D_t$, one should continue using the current
value of $D_t$, while in the set $\{(x,y): x\le -y B\}$ switch to the
opposite one. In what follows, we will call these sets the
\emph{continuation set} and the \emph{switching set}, respectively.

Next we need to find $B$. Introduce the value function $U(x,y)$
(cf.~\eqref{Ustar}; it turns out to be the same function $U(x,y)$ which
appears in the proof):
\[
U(x,y) = \inf_{D} \E_x  \biggl(\frac{c_1}{2} \int_{0}^\infty
(1-M_tD_t) dt + c_2\I(D_0\neq y) + c_2 \sum_{t> 0} \I(D_{t-} \neq
D_t) \biggr),
\]
where the infimum is taken over all c\`adl\`ag processes $D_t$ which are
adapted to the filtration generated by $M_t$, take on values $\pm1$, and
have a finite number of jumps. In the switching set, we have
\[
U(x,y) = U(x,-y)+c_2. 
\]
From the general theory (see Chapter~III in~\cite{PS06}), we can expect that
the value function $U(x,y)$ in the continuation set solves the ODE
\[
LU(x,y) = -\frac{c_1}{2}(1-xy).
\]
Its general solution can be found explicitly:
\[
U_{\text{gen}}(x,1) =  \frac{c_1(1-x)}{4\mu^2} \ln\frac{1+x}{1-x} + K_1x + K_2,
\]
where $K_1$ and $K_2$ are constants. Since we have $U(1,1)= 0$ (if $x=1$,
then $M_t=1$ for all $t\ge 0$ and the optimal process $D$ is $D_t\equiv 1$),
we get $K_2 = -K_1$. To find $K_1$ and $B$, we can employ the continuous fit
and smooth fit conditions, also known from the general theory, which state
that at the boundary of the continuation set, i.e.\ at the points $(x,y)$
with $x=-yB$, the value function satisfies the equations
\[
U(-B,1) = U(-B,-1) + c_2,\qquad U'(-B,1) = U'(-B,-1)
\]
(here $x=-B$, $y=1$; the pair $x=B$, $y=-1$ gives the same equations due to
the symmetry of the problem). Solving these equations gives formulas
\eqref{U-eq-1}--\eqref{U-eq-3} for $U(x,y)$.

To find the function $V(x)$ we use a similar approach. From the
representation as a standard optimal stopping problem \eqref{Vstar-3}, we
can expect that the optimal stopping time should be the first exit time of
the process $M_t$ from some continuation set. Taking into account the
original formulation of the problem as a sequential test, it is natural to
assume that the initial decision should be made at a moment when the
posterior mean becomes close to 1 or $-1$, i.e.\ the continuation set for
$V(x)$ should be an interval $(-A,\,A)$. As follows from the general theory,
$V(x)$ in the continuation set satisfies the ODE
\[
L V(x) = -c_0,
\]
which has the general solution 
\[
V_{\text{gen}}(x) = \frac{c_0x}{2\mu^2} \ln\frac{1-x}{1+x} + K_3 x + K_4.
\]
Due to the symmetry of the problem, we have $V(x) = V(-x)$, so $K_3 = 0$.
Then the constants $A$ and $K_4$ can be found from the continuous fit and
smooth fit conditions at $x=A$:
\[
V(A) = U(A,1), \qquad V'(A) = U'(A,1).
\]
These equations give the function $V(x)$ defined in
\eqref{V-eq-1}--\eqref{V-eq-2}, with $K_4=K$ from~\eqref{K}.
\end{remark}

\phantomsection
\addcontentsline{toc}{section}{References}%
\bibliographystyle{abbrv}
\bibliography{sequential-test-switching}

\end{document}